\magnification=1200

\overfullrule=0pt

\font\twelverm=cmr12
\font\twelvebf=cmbx12
\font\twelveit=cmti12
\def\bigfont{		
	\let\rm=\twelverm
	\let\bf=\twelvebf
	\let\it=\twelveit
	\rm}

\let\eu=\teneufm
\input amssym.def
\input amssym
\vsize=8truein
\hsize=6truein
\hoffset=18truept

\catcode`\@=11
\newskip\ttglue
\font\ninerm=cmr9

\font\sixrm=cmr6

\font\ninei=cmmi9
\font\eighti=cmmi8
\font\sixi=cmmi6
\skewchar\ninei='177 \skewchar\eighti='177 \skewchar\sixi='177

\font\ninesy=cmsy9
\font\eightsy=cmsy8
\font\sixsy=cmsy6
\skewchar\ninesy='60 \skewchar\eightsy='60 \skewchar\sixsy='60

\font\ninebf=cmbx9

\font\sixbf=cmbx6

\font\ninett=cmtt9
\font\eighttt=cmtt8

\hyphenchar\tentt=-1 
\hyphenchar\ninett=-1
\hyphenchar\eighttt=-1

\font\ninesl=cmsl9

\font\nineit=cmti9

\def\ninepoint{\def\rm{\fam0\ninerm}%
  \textfont0=\ninerm \scriptfont0=\sixrm \scriptscriptfont0=\fiverm
  \textfont1=\ninei \scriptfont1=\sixi \scriptscriptfont1=\fivei
  \textfont2=\ninesy \scriptfont2=\sixsy \scriptscriptfont2=\fivesy
  \textfont3=\tenex \scriptfont3=\tenex \scriptscriptfont3=\tenex
  \def\it{\fam\itfam\nineit}%
  \textfont\itfam=\nineit
  \def\sl{\fam\slfam\ninesl}%
  \textfont\slfam=\ninesl
  \def\bf{\fam\bffam\ninebf}%
  \textfont\bffam=\ninebf \scriptfont\bffam=\sixbf
   \scriptscriptfont\bffam=\fivebf
  \def\tt{\fam\ttfam\ninett}%
  \textfont\ttfam=\ninett
  \tt \ttglue=.5em plus.25em minus.15em
  \normalbaselineskip=11pt
  \def\MF{{\manual hijk}\-{\manual lmnj}}%
  \let\sc=\sevenrm
  \let\big=\ninebig
  \setbox\strutbox=\hbox{\vrule height8pt depth3pt width\z@}%
  \normalbaselines\rm}

\font\smaller=cmr8
\font\ninecsc=cmcsc9
\font\tencsc=cmcsc10

\def\beginsection#1{\bigskip\medskip
\centerline{\tencsc #1}\bigskip}

\catcode`\@=12

\headline={\ifnum\pageno>1\smaller\ifodd\pageno\hfill 
NEVANLINNA THEORY AND RATIONAL POINTS
\hfill\the\pageno \else
\the\pageno\hfill\uppercase{Junjiro Noguchi}\hfill\fi\else\hss\fi}

\footline={\hss}

\centerline{\bf NEVANLINNA THEORY AND RATIONAL POINTS}
\footnote{}{Research at MSRI supported in part by NSF grant \#DMS 9022140.}
\bigskip\medskip
\centerline{\bf Junjiro Noguchi}

\midinsert
\narrower\narrower
{\ninecsc Abstract}
\baselineskip=10pt
{\ninepoint 
S. Lang [L] conjectured in 1974 that a hyperbolic algebraic variety defined
over a number field has only finitely many rational points,
and its analogue over function fields.
We discuss the Nevanlinna-Cartan theory over function fields of
arbitrary dimension and apply it for Diophantine property of hyperbolic
projective hypersurfaces (homogeneous Diophantine equations)
constructed by Masuda-Noguchi [MN].
We also deal with the finiteness property of $S$-units points of
those Diophantine equations over number fields.

}
\endinsert

\beginsection{Introduction}

S. Lang [L] conjectured in 1974 that a hyperbolic algebraic variety defined
over a number field has only finitely many rational points,
and its analogue over function fields.
For subvarieties of Abelian varieties the function field analogue was dealt
with by M. Raynaud [R], and lately G. Faltings [F] proved this conjecture for
subvarieties of Abelian varieties over number fields.
On the other hand, the author [Nog5] proved the function field analogue
in general case (cf.\ also [Nog1], [Nog2]).
See Y. Imayoshi-H. Shiga [IS], M. Zaidenberg [Z], and
M. Suzuki [Su1], [Su2] for non-compact versions of this result.

In the case of curves (Fermat, Catalan, Thue equations, etc.)
defined over function fields,
R.C. Mason [Ma1], J. Silverman [Si] and J. Mueller [Mu] obtained
similar or more precise finiteness properties
by making use of a different method which relies on the function field
analogue of ``$abc$-conjecture" of Masser-Oesterl\'e.
The function field analogue of ``$abc$-conjecture"
was proved in more general form by R.C. Mason [Ma1], [Ma2],
J. Voloch [Vo] and W. Brownawell-D. Masser [BrM].
They actually proved a version of
``$abc$-conjecture" in several variables (to say, $abc\cdots$-conjecture),
which is nothing but a special case of
Nevanlinna-Cartan's second main theorem with truncated counting functions
applied to algebraic case (see [C], (3), and \S2).

Here we discuss the Nevanlinna-Cartan theory over function fields of
arbitrary dimension (cf.\ J. Wang [W1], [W2] and [W3] for related
results), and apply it for Diophantine property of hyperbolic
projective hypersurfaces (homogeneous Diophantine equations)
constructed by [MN].
We also deal with the finiteness property of $S$-units points of
those Diophantine equations over number fields.

{\it Acknowledgement.}
The present manuscript is based on a talk of the author for
Workshop in Nevanlinna Theory and Diophantine Approximation, January 22-26,
1996, Mathematical Sciences Research Institute,
University of California, Berkeley.
The author is grateful to the institute for the very active
and stimulating circumstance, and for the hospitality.
His thanks are also due to the co-organizers of the workshop,
Professor Pit-Mann Wong and P. Vojta.

\beginsection{1. Function field case}

We deal with the Nevanlinna theory over function fields of
an arbitrary dimension by making use of the method of Cartan [C] and
Nochka [N], which affirmatively proved Cartan's conjecture
(cf.\ also [Ch] and [Fu2]).
Let $k$ be an algebraically closed field of characteristic 0,
let $R$ be a smooth projective algebraic variety of dimension $N$ over $k$,
and let $K$ denote the rational function field of $R$.
To use analytic definitions as well, we set $k={\bf C}$.

We fix a class $\omega$ of an ample line bundle (a Hodge metric form) on $R$.
Let $D$ be a divisor on $R$, and define the counting function of $D$
with respect to $\omega$ by
$$
N(D;\omega)=\int_D \omega^{N-1}.
$$
If $\sigma$ is a meromorphic section of some line bundle and $(\sigma)$ denotes
the divisor determined by $\sigma$, then we write
$$
N(\sigma;\omega)=N((\sigma);\omega).
$$
Let $a_i \in K^*$, $1 \leqq i \leqq m$, and let $((a_j))$ denote the least
common multiple of the polar divisors of $a_i$, $1 \leqq i \leqq m$.
We define the height function of $(a_i)_{1 \leqq i \leqq m}$, by
$$
{\rm ht}((a_j);\omega)=N(((a_i)); \omega).
$$
From $a_i \in K^*$, $1 \leqq i \leqq m$, we obtain a rational mapping
$f=[1; a_1, \ldots , a_m]:R \rightarrow {\bf P}^m({\bf C})$, and define
the order function of $f$ by
$$
T(f;\omega)=\int_R f^*\Omega \wedge \omega^{N-1},
$$
where $\Omega$ denotes the Fubini-Study K\"ahler form on ${\bf P}^m({\bf C})$.
The first main theorem is nothing but the Poincar\'e residue theorem:

{\bf (1.1) Theorem.}
{\it Let the notation be as above. Let $L$ be the line bundle determined by
$(a_i)$, and let $\sigma \in \Gamma(R,L)$ be a holomorphic section.
Then
$$
T(f;\omega)={\rm ht}((a_j);\omega)=N(\sigma;\omega),\qquad 0 \leqq j \leqq m.
$$
}

Next we have the
``{\it second main theorem with truncated counting functions}'':

{\bf (1.2) Theorem.}
{\it
Let $f=[\sigma_0,\ldots,\sigma_m]: R \rightarrow {\bf P}^m ({\bf C})$
be a reduced representaion of a rational mapping given
by holomorphic sections $\sigma_j$ of a line bundle.
Let $H_1,\ldots, H_q$ $(q \geqq m+1)$ be linear forms on ${\bf P}^m ({\bf C})$
in general position such that $f^*H_j \not\equiv 0$ for any $j$.
Let $r$ denote the rank of $df$ at general point, and let $l$ denote
the dimension of the smallest linear subspace of ${\bf P}^m ({\bf C})$
containing $f(R)$.
Then
$$
\eqalign{
(q-l-1)T(f;\omega) \leqq &\sum_{i=1}^q N_{l-r+1}(f^*H_i;\omega)\cr
&+\left\{{{(l-r+1)(l-r+2)}\over 2} +r-1\right\}
N(J; \omega).\cr}
\leqno(1.3)
$$
}

Here $N_{l-r+1}(H_i(x_j);\omega)$ is the truncated counting function
of zeros of $H_i(x_j)$, and $J$ is a divisor used to define
the {\it generalized Wronskian bundle} $W=[pJ]\otimes L^{l+1}$,
and is independent of $f$.
Cf.\ J. T.-Y. Wang [W1], [W2], and [W3] for related results.

Theorem (1.2) plays an important role.
Because of the truncated counting functions $N_{l-r+1}(H_i(x_j);\omega)$,
we may consider inequality (1.3) as a version of
``$abc$-conjecture'' in several variables
over function fields of an arbitrary dimension.
(cf.\ [Ma2], [Vo], and [BrM], Corollary I).
We derive a ``ramification theorem'' over function fields (Corollary (2.16)),
and ``Borel's theorem'' over functions fields:

{\bf (1.4) Theorem.}
{\it
Let $x_j \in K^*, 1 \leqq j \leqq s$, satisfy
$$
a_1x_1^d + \cdots +a_sx_s^d=0\qquad(s \geqq 2).
$$
Assume
$$
d > s(s-2)+(s-1)^2{\rm ht}((a_1,\ldots, a_s); \omega)
+{{(s-1)(s-2)}\over 2}N(J; \omega).
\leqno(1.5)
$$
Then there is a disjoint decomposition
$\{1,\ldots,s\}=\cup_{\nu=1}^l I_\nu$ of indices such that
\parindent=30pt

\item{\rm (i)}
$|I_\nu| \geqq 2$ for all $\nu$;
\item{\rm (ii)}
for arbitrary two indices $j, k \in I_\nu$, the ratio $x_j/x_k$ is a constant.
\item{\rm (iii)}
$\sum_{j \in I_\nu} a_j x_j^d=0$ for all $\nu$.
}
\parindent=10pt

We use the above results to study the rational points
of a Diophantine equation.
Let $X \subset {\bf P}_K^{n-1}$
be a hypersurface defined over $K$ by equation
$$
a_1M_1^d(z_1,\ldots,z_n)+\cdots+a_sM_s^d(z_1,\ldots,z_n)=0,
\leqno(1.6)
$$
where $a_j \in K^*=K\setminus \{0\}$ and $d \in {\bf Z},>0$.
We set
$$
\eqalign{
Y({\rm P})=\Bigl\lbrace(u_1,\ldots,u_n) \in {\bf P}^{n-1}(k);&
\sum_j a_j M_j^d(z_1,\ldots,z_n)M_j^d(u_1,\ldots,u_n)=0,\cr
&\hbox{ and } u_j=0 \hbox{ if } z_j=0 \Bigr\rbrace .\cr}
$$
Then $Y({\rm P})$ is a projective variety defined over $k$.
Moreover, we set
$$
{\cal R}({\rm P})=\left\{(z_1 u_1,\ldots,z_n u_n) \in {\bf P}^{n-1}(K);
(u_1,\ldots,u_n) \in Y({\rm P})\right\} \subset X(K).
$$

{\bf (1.7) Main Theorem.}
{\it
Let the notation be as above.
Assume that $\{M_j(z_1,\ldots,$ $z_n)\}_{j=1}^s$ is $n$-admissible ([MN]).

{\rm (i)}
Assume that $d > s(s-2)$.
Then the heights ${\rm ht}((z_i);\omega)$ of points of $X(K)$ are
bounded.

{\rm (ii)}
Assume (4) for $d$.
Then all points of $X(K)$ are defined over $k$; that is,
all points $x=(x_1,\ldots,x_n)\in X(K)$ are represented by $x_j \in k$.

{\rm (iii)}
Assume that
$$
d >s!(s!-2)+{{(s!-1)(s!-2)}\over 2}N(J; \omega).
$$
Then there are only finitely many rational points ${\rm P}_\mu \in X(K)$,
$\mu=1,\ldots,\mu_0 \; (< \infty)$ such that
$$
X(K)=\bigcup_{\mu=1}^{\mu_0} {\cal R}({\rm P_\mu}).
$$
}

In the proof of (iii) we use the result of E. Bombieri-J. Mueller [BM]
in a form generalized by Theorem (1.2).

\beginsection{2. Number field case--$S$-unit points}

In this section we deal with $X$ defined by (1.6) over a number field $F$.
We say that a point $(z_i) \in {\bf P}^{n-1}(F)$ is $S$-{\it units point}
if $z_i \in {\cal O}_S^*$ or $z_i=0$, denote
by $X({\cal O}_S^*)$ the set of $S$-units point of $X \subset {\bf P}^{n-1}_F$.
By making use the same idea as in the proof of the Main Theorem (1.7)
we then apply Borel's Theorem for $S$-units to prove

{\bf (2.1) Theorem.}
{\it
Let $\{M_j(z_1,\ldots,z_n)\}_{j=1}^s$ be an $n$-admissible set
of monomials, and let $X$ be defined by (1.6) with an arbitrary $d \geq 1$.
Then $X({\cal O}_S^*)$ is a finite set.
}

{\it Remark.}
(i) Mahler ([M], p.\ 724, Folgerung 2) proved Theorem (2.1)
in the case of $n=s=3$.

(ii) Let $[F;{\bf Q}]$ denote the extension degree and let $|S|$ be the
cardinality of $S$.
By making use of the bound obtained by Gy\"ory [G1, G2], we have
$$
|X({\cal O}^*_S)| \leqq (2^n -1) (2^{35} s^2)^{(s-1)^3 |S|}.
$$

(iii)(ASMT-conjecture) We give an analogue of Theorem (2.1) over number fields,
which is an extension of Schmidt's linear subspace theorem.
This may be called ``{\it Arithmetic Second Main Theorem-conjecture}\/''
by its nature.
Let $F$ and $S$ be as above.
We define
the truncated counting function $N_\lambda(H(x_i))$ over the
places of $F$ outside $S$.
Let $H_j$, $1 \leqq j \leqq q$, be linear forms in general
position on ${\bf P}^m_F$.
Then for an arbitrary $\epsilon >0$ there are finitely many hyperplanes
$E_\nu$ such that for $(x_i) \in {\bf P}^m(F)\setminus \cup E_\nu$
$$
(q-m-1-\epsilon){\rm ht}((x_i)) \leqq N_m(H(x_i)).
$$

(iv) The above ASMT-conjecture implies the finiteness of $X(F)$.
In this respect the works of Khoai-Tu [KT] and Sarnak-Wang [SW] are
of interest; especially Sarnak-Wang [SW] showed that some hyperbolic
smooth hypersurfaces of ${\bf P}^4$ and ${\bf P}^5$ constructed by [MN]
defined over $\bf Q$ has infinitely many rational points over
the $p$-adic number filed ${\bf Q}_p$
for every prime $p$, and that its Brauer-Manin group known as an
obstruction for the Hasse principle is vanishing.

(v) See [MN] for a number of examples.
\parindent=25pt

\beginsection{\centerline{References}}

\ninepoint

\item{[BM]}
E. Bombieri and J. Mueller,
The generalized Fermat equation in function fields,
J. Number Theory
{\bf 39} (1991), 339-350.
\item{[BrM]}
W.D. Brownawell and D.M. Masser,
Vanishing sums in function fields,
Math.\ Proc.\ Camb.\ Phil.\ Soc.\
{\bf 100} (1986), 427-434.
\item{[C]}
H. Cartan,
Sur les z\'eros des combinaisons lin\'eaires de $p$ fonctions holomorphes donn\'ees,
Mathematica {\bf 7} (1933), 5-31.
\item{[Ch]}
W. Chen,
Defect relations for degenerate meromorphic maps,
Trans.\ Amer.\ Math.\ Soc.\ {\bf 319} (1990), 499-415.
\item{[F]}
G. Faltings,
Diophantine approximation on Abelian varieties,
Ann.\ Math.\ {\bf 133} (1991), 549-576.
\item{[Fu1]}
H. Fujimoto,
Non-integrated defect relation for meromorphic maps of complete K\"ahler manifolds into ${\bf P}^{N_1}({\bf C}) \times \cdots \times {\bf P}^{N_k}({\bf C})$,
Japan.\ J. Math.\ {\bf 11} (1985), 233-264.
\item{[Fu2]}
---\hglue-1pt---,
Value Distribution Theory of the Gauss Map of Minimal Surfaces in ${\bf R}^m$,
Aspect Math.\ {\bf E21}, Vieweg, Braunschweig, 1993.
\item{[G1]}
K. Gy\"ory, On the numbers of families of solutions of systems of decomposable
form equations,
Publ.\ Math.\ Debrecen {\bf 42} (1993), 65-101.
\item{[G2]}
---\hglue-1pt---,
Letter, October 1994.
\item{[IS]}
Y. Imayoshi and H. Shiga,
A finiteness theorem for holomorphic families of Riemann surfaces,
In: D.\ Drasin (ed.)
Holomorphic Functions and Moduli vol.\ II,
Springer, New York-Berlin-Heidelberg-London-Paris-Tokyo, 1988.
\item{[KT]}
H. H. Khoai and M. V. Tu,
$p$-adic Nevanlinna-Cartan theorem, Internat.\ J. Math.\ {\bf 6} (1995),
719-731.
\item{[Ko1]}
S. Kobayashi,
Hyperbolic Manifolds and Holomorphic Mappings,
Pure and Appl.\ Math.\ vol.\ 2, Marcel Dekker, New York, 1970.
\item{[Ko2]}
---\hglue-1pt---,
Intrinsic distances, measures, and geometric function theory.
Bull.\ Amer.\ Math.\ Soc.\ {\bf 82} (1976), 357-416.
\item{[L1]}
S. Lang,
Higher dimensional Diophantine problems,
Bull.\ Amer.\ Math.\ Soc.\ {\bf 80} (1974), 779-787.
\item{[L2]}
---\hglue-1pt---,
Hyperbolic and Diophantine analysis.
Amer.\ Math.\ Soc.\ {\bf 14} (1986), 159-205.
\item{[L3]}
---\hglue-1pt---,
Introduction to Complex Hyperbolic Spaces,
\hyphenation{-}
Springer-Verlag, New York-Berlin-Heidelberg, 1987.
\item{[L4]}
---\hglue-1pt---,
Number Theory III, Encycl.\ Math.\ Sci.\ vol.\ {\bf 60},
Springer-Verlag,
Berlin-Heidel-berg-New York-London-Paris-Tokyo-Hong Kong-Barcelona, 1991.
\item{[La]}
M. Laurent, Equations diophantinnes exponentielles,
Invent.\ Math.\ {\bf 78} (1984), 299-327.
\item{[M]}
K. Mahler,
Zur Approximation algebraischer Zahlen.\ I.\ (\"Uber den gr\"o\ss ten Primteiler bin\"arer Formen.),
Math.\ Ann.\ {\bf 107} (1983), 691-730.
\item{[Ma1]}
R.C. Mason,
Diophantine Equations over Function Fields,
London Math.\ Soc.\ Lecture Notes vol.\ {\bf 96},
Cambridge University Press,  Cambridge, 1984.
\item{[Ma2]}
---\hglue-1pt---,
Norm form equations.\ I,
J. Number Theory {\bf 22} (1933), 190-207.
\item{[Mu]}
J. Mueller,
Binomial Thue's equation over function fields,
Compo.\ Math.\ {\bf 73} (1990), 189-197.
\item{[MN]}
K. Masuda and J. Noguchi,
A construction of hyperbolic hypersurfaces of ${\bf P}^n ({\bf C})$,
Math.\ Ann.\ {\bf 304} (1996), 339-362.

\item{[N]}
A.M. Nadel,
The nonexistence of certain level structures on abelian varieties over complex function fields,
Ann.\ Math.\ {\bf 129} (1989), 161-178.
\item{[No]}
E.I. Nochka, On the theory of meromorphic functions, Sov.\ Math.\ Dokl.\
{\bf 27} (1983), 377-381.
\item{[Nog1]}
J. Noguchi,
A higher dimensional analogue of Mordell's conjecture over function fields,
Math.\ Ann.\ {\bf 258} (1981), 207-212.
\item{[Nog2]}
---\hglue-1pt---,
Hyperbolic fibre spaces and Mordell's conjecture over function fields,
Publ.\ RIMS, Kyoto Univ.\ {\bf 21} (1985), 27-46.
\item{[Nog3]}
---\hglue-1pt---,
Hyperbolic manifolds and Diophantine geometry,
Sugaku Expositions {\bf 4} (1991), pp.\ 63-81, Amer.\ Math.\ Soc.,
Providence, Rhode Island, 1991.
\item{[Nog4]}
---\hglue-1pt---,
Moduli space of Abelian varieties with level structure over function fields,
International J. Math.\ {\bf 2} (1991), 183-194.
\item{[Nog5]}
---\hglue-1pt---,
Meromorphic mappings into compact hyperbolic complex spaces and geometric Diophantine problem,
Internat.\ J. Math.\ {\bf 3} (1992), 277-289.
\item{[Nog6]}
---\hglue-1pt---,
An example of a hyperbolic fiber space without hyperbolic embedding into compactification,
Proc.\ Osaka International Conference, Osaka 1990, Complex Geometry
(G.\ Komatsu and Y.\ Sakane, eds.)
Lecture Notes in Pure and Appl.\ Math.\ vol.\ {\bf 143}, pp.\ 157-160,
Marcel Dekker, New York-Basel-Hong Kong, 1993.
\item{[Nog7]}
---\hglue-1pt---,
Some topics in Nevanlinna theory, hyperbolic manifolds and Diophantine geometry,
Geometry and Analysis on Complex Manifold, Festschrift for A. Kobayashi 60th Birthday,
pp. 140--156, World Scientific, Singapore, 1994.
\item{[NO]}
---\hglue-1pt--- and T. Ochiai,
Geometric Function Theory in Several Complex Variables,
Transl.\ Math.\ Mono.\ {\bf 80}, Amer.\ Math.\ Soc.,
Providence, Rhode Island, 1990.
\item{[R]}
M. Raynaud,
Around the Mordell conjecture for function fields and a conjecture of S. Lang,
Algebraic Geometry (M. Raynaud and T. Shioda, eds.),
Lecture Notes in Math.\ vol.\ 1016, pp.\ 1-19,
Springer, Berlin-New York, 1983.
\item{[SW]}
P. Sarnak and L. Wang,
Some hypersurfaces in ${\bf P}^4$ and the Hasse-principle,
C.R. Acad.\ Sci.\ Paris, S\'er.\ I {\bf 321} (1995), 319-322.
\item{[S]}
H.P. Schlickewei,
The {\eu p}-adic Thue-Siegel-Roth-Schmidt Theorem,
Archiv
\break
Math. {\bf 29} (1977), 267-270.
\item{[Sc]}
W. Schmidt,
Diophantine Approximation,
Lecture Notes in Math.\ vol.\ {\bf 785},
Springer, Berlin-Heidelberg-New York, 980.
\item{[Si]}
J.H. Silverman,
The Catalan equation over function fields,
Trans.\ Amer.\ Math.\ Soc.\ {\bf 273} (1982), 201-205.
\item{[St1]}
W. Stoll,
Die beiden Haupts\"atze der Wertverteilungstheorie bei Funktionen mehrerer komplexer Ver\"anderlichen (I); (II),
J. Acta.\ Math.\ {\bf 90}, (1953),1-115; Acta Math.\b {\bf 92} (1954), 55-169.
\item{[St2]}
---\hglue-1pt---,
Value Distribution of Holomorphic Maps into Compact Complex Manifolds,
Lecture Notes in Math.\ {\bf 135}, Springer-Verlag,
Berlin-Heidelberg-New York, 1970.
\item{[St3]}
---\hglue-1pt---,
Value Distribution on Parabolic Spaces,
Lecture Notes in Math.\ {\bf 600},
Springer-Verlag, Berlin-Heidelberg-New York, 1977.
\item{[Su1]}
M. Suzuki,
Moduli spaces of holomorphic mappings into hyperbolically embedded complex spaces and holomorphic fibre spaces,
to appear in J. Math.\ Soc.\ Japan.
\item{[Su2]}
---\hglue-1pt---,
Mordell property of hyperbolic fiber spaces with noncompact fibers,
to appear in T\^ohoku J. Math.
\item{[V]}
P. Vojta,
Diophantine Approximations and Value Distribution Theory,
Lecture Notes in Math.\ vol.\ {\bf 1239},
Springer, Berlin-Heidelberg-New York, 1987.
\item{[Vo]}
J.F. Voloch,
Diagonal equations over function fields,
Bol.\ Soc.\ Brasil.\ Math.\ {\bf 16} (1985), 29-39.
\item{[W1]}
J. T.-Y. Wang,
Diophantine equations over function fields,
Thesis, Notre Dame,  1994.
\item{[W2]}
 ---\hglue-1pt---,
The truncated second main theorem of function fields,
to appear in J. Number Theory.
\item{[W3]}
 ---\hglue-1pt---,
Effective Roth theorem of function fields,
The Rocky Mountain J. Math.

\item{[Z]}
M.G. Zaidenberg,
A function-field analog of the Mordell conjecture: a noncompact version,
Math.\ USSR Izvestiya {\bf 35} (1990), 61-81.

\bigskip\rm
\noindent{\ninecsc
Department of Mathematics, Tokyo Institute of Technology, Ohokayama, Meguro,
Tokyo 152, Japan ({\tt noguchi@math.titech.ac.jp}).}

\medskip{\it \noindent Current address:}
Mathematical Sciences Research Institute, University of California, Berkeley,
1000 Centennial Drive, Berkeley, California 94720, USA
({\tt noguchi@\allowbreak msri.org}).

\bye